\documentclass[11pt]{article}

\usepackage{latexsym, amsmath}
\newtheorem{theorem}{Theorem}
\newtheorem{lemma}{Lemma}
\newtheorem{corollary}{Corollary}

\newtheorem{conjecture}{Conjecture}
\newtheorem{proposition}{Proposition}
\newtheorem{algorithm}{Algorithm}
\newenvironment{proof}{{\it Proof:\/}}{\hfill $\Box$\\ }
 
\newcommand{\Z}{\mbox{\bf Z}}
\newcommand{\F}{\mbox{\bf  F}}

\newcommand{\floor}[1]{\lfloor #1 \rfloor}
\newcommand{\ceil}[1]{\lceil #1 \rceil}

\title{On the Bounded Sum-of-digits Discrete Logarithm Problem
in Kummer and Artin-Schreier Extensions}

\date{}

\author{Qi Cheng\thanks{School of Computer Science,
the University of Oklahoma,
Norman, OK 73019, USA.
Email: {\tt qcheng@cs.ou.edu.}
This research is partially supported by NSF Career
Award CCR-0237845.
}}

\begin{document}

\maketitle

\begin{abstract}

In this paper, we study the discrete logarithm problem in the finite
fields $\F_{q^n}$ where $n|q-1$.  The field is called a Kummer field
or a Kummer extension of $\F_q$. It plays an important role in
improving the AKS primality proving algorithm.  It is known that we
can efficiently construct an element $g$ with order greater than $2^n$
in the fields.  Let $S_q(\bullet)$ be the function from integers to
the sum of digits in their $q$-ary expansions.  We present an
algorithm that given $g^e$ ($ 0\leq e < q^n $ ) finds $e$ in random
polynomial time, provided that $S_q (e) < n$.  We
then show that the problem is solvable in random polynomial time for
most of the exponent $e$ with $S_q (e) < 1.32 n $.  The main tool for
the latter result is the Guruswami-Sudan list decoding algorithm.
Built on these results, we prove that in the field $\F_{q^{q-1}}$, the
bounded sum-of-digits discrete logarithm with respect to $g$
can be computed in random time $O(f(w)  \log^4 (q^{q-1}))$, 
where $f$ is a subexponential function and $w$ is
the bound on the $q$-ary sum-of-digits of the exponent.  Hence the
problem is fixed parameter tractable.  These results are shown to be
extendible to Artin-Schreier extension $\F_{p^p}$ where $p$ is a
prime.  Since every finite field has an extension of reasonable degree
which is a Kummer field, our result reveals an unexpected property of
the discrete logarithm problem, namely, the bounded
sum-of-digits discrete logarithm problem in any given finite field
becomes polynomial time solvable in 
certain low degree extensions.
\end{abstract}

\section{Introduction and Motivations}

Most of practical public key cryptosystems base their security on the
hardness of solving the integer factorization problem or the discrete
logarithm problem in finite fields. Both of the problems admit
subexponential algorithms, thus we have to use long parameters, which
make the encryption/decryption costly if the parameters are randomly
chosen.  Parameters of low Hamming weight, or more generally, of
small sum-of-digits, offer some remedy. Using them speeds up the
system while seems to keep the security intact.  In particular, in the
cryptosystem based on the discrete logarithm problem in finite fields
of small characteristic, using small sum-of-digits exponents is very
attractive, due to the existence of normal bases \cite{AgnewMu91}.  It
is proposed and implemented for smart cards and mobile devices, where
the computing power is severely limited.  Although attacks exploring
the specialty were proposed \cite{Stinson02}, none of them have
polynomial time complexity.

Let $\F_{q^n}$ be a finite field. 
For $\beta\in \F_{q^n}$, if $\beta, \beta^q, \beta^{q^2}, \cdots,
\beta^{q^{n-1}} $ form a linear basis of $\F_{q^n}$ over $\F_q$, we
call them {\em a normal basis}. It is known that a normal basis exists
for every pair of prime power $q$ and a positive integer $n$
\cite[Page 29]{Gao93}. Every element $\alpha $ in $\F_{q^n}$ can be
represented as
$$ \alpha = a_0 \beta + a_1 \beta^q + \cdots + a_{n-1} \beta^{q^{n-1}} $$
where $a_i \in \F_q $ for $0\leq i\leq n-1$.
The power of $q$ is a linear operation,
$$ \alpha^q =   a_0 \beta^q + \cdots + a_{n-2}
\beta^{q^{n-1}} + a_{n-1} \beta. $$ 
Hence to compute the $q$-th power, we only need to
shift the digits, which can be done very fast, possibly
on the hardware level. Now suppose we want to
compute $ \alpha^e $ where the $q$-ary expansion of $e$ is
\begin{equation} \label{expansionofe}
e= e_0 + e_1 q + e_2 q^2 + \cdots + e_{n-1} q^{n-1}\ \ \ \ (0\leq e_i <
q\ {\rm for}\ \ 0\leq i\leq n-1). 
\end{equation}
The sum-of-digits of $e$ in the $q$-ary expansion is defined as $S_q
(e) = \sum_{i=0}^{n-1} e_i$. When $q=2$, the sum-of-digits becomes the
famous Hamming weight.  To compute  $\alpha^e$, we
only need to do shiftings and at most $S_q (e)$ number of
multiplications. Furthermore, the exponentiation algorithm can be
parallelized,
which is a property not enjoyed by the large
characteristic fields. For details, see \cite{Gathen91}.

\subsection{Related work}

The discrete logarithm problem in finite field $\F_{q^n}$, is to
compute an integer $e$ such that $g' = g^e$, given a generator $g$ of
a subgroup of $\F_{q^n}^*$ and $g'$ in the subgroup.  The general
purpose algorithms to solve the discrete logarithm problem are the number
field sieve and the function field sieve (for a survey see
\cite{Odlyzko00}). They have time complexity
 $$ {\rm exp}(c (\log q^n)^{1/3} (\log \log q^n)^{2/3} ) $$ for some
constant $c$, when $q$ is small, or $n$ is small.  

Suppose we want to
compute the discrete logarithm of $g^e$ with respect to base $g$ in the
finite field $\F_{q^n}$.  If we know that the Hamming weight of $e$ is
equal to $w$, there is an algorithm proposed by Coppersmith, which
works well if $w$ is very small. It is a clever adaption of the
baby-step giant-step idea, and runs in random time $O(\sqrt{w} {
\lfloor \log q^n /2 \rfloor \choose \lfloor w/2 \rfloor } )$.  It is
proved in \cite{Stinson02} that the average-case complexity achieves
only a constant factor speed-up over the worst case.  It is not clear
how his idea can be generalized when the exponent has small
sum-of-digits in the base $q > 2$. However, we can consider the very
special case where $e_i \in \{0,1\}$ for $0\leq i\leq n-1$ and
$\sum_{0\leq i\leq n-1} e_i = \lfloor {n \over 2} \rfloor$, recall
that $e_i$'s are the digits of $e$ in the $q$-ary expansion. It can be
verified that the Coppersmith algorithm can be applied in this
case. The time complexity becomes $O( \sqrt{n} { \lfloor n/2 \rfloor
\choose \lfloor n/4\rfloor})$. If $q < n^{O(1)}$, it is much worse
than the time complexity of the function field sieve on a general
exponent.

If the $q$-ary sum-of-digits of the exponent is bounded by $w$,
is there an algorithm which runs in time $f(w) \log^c (q^n)$ and
solves the discrete logarithm problem in $\F_{q^n}$, for a arbitrary
function $f$  and a constant $c$?
A similar problem has been raised from the parametric point of view
by Fellows and Koblitz \cite{FellowsKo93}, where they consider the
prime finite fields and the bounded Hamming weight exponents.  
Their problem is listed among the most important open problems in the
theory of parameterized complexity \cite{DowneyFe99}.  From the above
discussions, it is certainly more relevant to cryptography to treat the
finite fields with small characteristic and exponents with bounded
sum-of-digits.

Unlike the case of the integer factorization, where a lot of special
purpose algorithms exist, the discrete logarithm problem is considered
more intractable in general.  As an example, one should not use a RSA
modulo of about 1000 bits with one prime factor of 160 bits. It would
be vulnerable to the elliptic curve factorization algorithm. However,
in the Digital Signature Standard, adopted by the U.S. government, the
finite field has cardinality about $2^{1024}$ or larger, while the
encryption/decryption is done in a subgroup of cardinality about
$2^{160}$.  As another example, one should search for a secret prime
as random as possible in RSA, while in the case of the discrete
logarithm problem, one may use a finite field of small characteristic,
hence the group of very special order.  It is believed that no
trapdoor can be placed in the group order, as long as it has a large
prime factor (see the panel report on this issue in the Proceeding of
Eurocrypt 1992).  In order to have an efficient algorithm to solve the
discrete logarithm, we need that every prime factor of the group order
is bounded by a polynomial function on the cardinality of the field.
Given the current state of analytic number theory, it is very hard, if
not impossible, to prove that there exists {\em infinite} many of
finite fields of even (or constant) characteristic, where the discrete
logarithm can be solved in polynomial time.

In summary, there are several common perceptions about the discrete
logarithm problem in finite fields:
\begin{enumerate}
\item As long as the group order has a big prime factor, the discrete
logarithm problem is hard. We may use exponents with small
sum-of-digits, since the discrete logarithm problem in that case seems to be
fixed parameter intractable. We gain advantage in speed 
by using bounded
sum-of-digits exponents, and at the same time keep 
the problem as infeasible as using the general exponents.
\item If computing discrete logarithm is difficult, it should be
difficult for any generator of the group.  The discrete logarithm
problem with respect to one generator can be reduced to the discrete
logarithm problem with respect to any generator. Even though in the small
sum-of-digits case, a reduction is not available, it is
not known that changing the generator of the group affects the
hardness of the discrete logarithm problem.
\end{enumerate}

\subsection{Our results}

In this paper, we show that those perceptions are problematic,
by studying the discrete logarithm problem in large multiplicative 
subgroups of the Kummer and Artin-Schreier extensions with a 
prescribed generator.
We prove that the bounded sum-of-digits discrete logarithm 
are easy in those groups. 
More precisely we prove constructively:
\begin{theorem}\label{main}
(Main) There exists a random algorithm to find the integer $e$
given $g$ and $g^e$ in $\F_{q^n}$ in time polynomial in $\log (q^n)$
under the conditions:
\begin{enumerate}
\item $n|q-1$;
\item $0\leq e< q^n$, and $S_q (e) \leq n$;
\item \label{conditiong} $g = \alpha+b $ where $\F_q (\alpha) =
\F_{q^n}$, $b\in \F_q^* $ and $\alpha^n \in \F_q$.
\end{enumerate}
Moreover, there does not exist an integer $e' \not= e$  
satisfying that $0\leq e'< q^n$, $S_q (e') \leq n$
and $g^{e'} = g^e$
\end{theorem}

A few comments are in order: 
\begin{itemize}
\item For a finite field $\F_{q^n}$, if $n|q-1$, then there exists $g
\in \F_{q^n}$ satisfying the condition in the theorem;  if
there exists $\alpha$ such that $\F_q (\alpha) =
\F_{q^n}$ and $\alpha^n \in \F_q$, then $n | q-1$.
\item As a comparison, Coppersmith's algorithm runs in exponential
time in the case where $e_i \in \{0,1\}$ for $0\leq i\leq n-1$, 
$S_q (e) = { n \over 2} $ and $q < n^{O(1)}$, while our algorithm runs in
polynomial time in that case. On the other hand, Coppersmith's
algorithm works for every finite field, while our algorithm works in
Kummer fields. Our result has an indirect affect on an arbitrary
finite field though, since every finite field has extensions of degree close
to a given number, which are Kummer fields. As an example, suppose we
want to find an extension of $\F_q$ with degree about $\log^2 q$.  We
first pick a random $n$ close to $\log q$ such that $(n,q) =1$. Let
$l$ be the order of $q$ in $\Z/n\Z$. The field $\F_{(q^l)^n}$ is a
Kummer extension of $\F_{q^l}$, and an extension of $\F_q$. According
to Theorem~\ref{main}, there is a polynomial time algorithm which
computes the discrete logarithm to some element $g$ in $\F_{q^{ln}}$
provided that the sum-of-digits of the exponent in the $q^l$-ary
expansion is less than $n$.  Hence our result reveals an unexpected
property of the discrete logarithm problem in finite fields: the
difficulty of bounded sum-of-digits discrete logarithm problem drops
dramatically if we move up to extensions.
\item Numerical evidences suggest that the order of $g$ is
close to the group order $q^n -1 $, if it does not equal to $q^n-1$.
However, it seems hard to prove it.  In fact, this is one of the
main obstacles in improving the efficiency of AKS-style primality testing
algorithm \cite{AgrawalKa02}.  We make the following conjecture.
\begin{conjecture} \label{orderconjecture}
Suppose that a finite field $\F_{q^n}$ and an element $g$ in the field 
satisfy the conditions in 
Theorem~\ref{main}. In addition, $ n \geq \log q$.
The order of $g$ is greater than $q^{n/c}$ for an absolute constant $c$.
\end{conjecture}
\item Even though we can not prove that the largest prime factor of
the order of $g$ is very big, it seems, as supported by numerical
evidences, that the order of $g$, which is a factor of $q^n -1 $
bigger than $2^n$, is rarely smooth.  For instance, in the
$\F_{2^{889}} = \F_{128^{127}}$, any $g$ generates the whole
group $\F_{2^{889}}^*$.  The order $2^{889} - 1$ contains a prime
factor of $749$ bits.  One should not attempt to apply the
Silver-Pohlig-Hellman algorithm here.
\end{itemize}

A natural question arises: can the restriction on the sum-of-digits in
Theorem~\ref{main} be relaxed?  Clearly if we can solve the problem
under condition $S_q (e) \leq (q-1) n$ in polynomial time, then the
discrete logarithm problem in subgroup generated by $g$ is broken. If
$g$ is a generator of $\F_{q^n}^*$, then the discrete logarithm
problem in $\F_{q^n}$ and any of its subfields to any base are broken.
We find a surprising relationship between the relaxed problem and the
list decoding problem. We are able to prove:

\begin{theorem}\label{listdecoding}
Suppose $e$ is chosen in random  from the set
$$ \{ 0\leq e < q^n -1 | S_q (e) < 1.32 n \}. $$ There exists an algorithm
given $g$ and $g^e$ in $\F_{q^n}$, to find $e$ in time polynomial in
$\log (q^n)$, with probability greater than $1 - c^{-n}$ for some
constant $c$ greater than $1$, under the conditions:
\begin{enumerate}
\item $n|q-1$;
\item $g = \alpha+b $ where $\F_q (\alpha) = \F_{q^n}$, $b\in \F_q^* $
and $\alpha^n \in \F_q$.
\end{enumerate}
\end{theorem}

We also prove a parameterized complexity result concerning
the bounded sum-of-digits discrete logarithm.

\begin{theorem}\label{parameter}
There exists an element $g$ of order greater than $2^q$ in 
$\F_{q^{q-1}}^*$, such that the discrete logarithm problem 
with respect to the generator $g$
can be solved in time $f(w) \log^4 (q^{q-1})$,
where $f$ is a subexponential function and
$  w $ is the bound of the sum-of-digits of the 
exponent in $q$-ary expansion.
\end{theorem}

This answers an important open question in
parameterized complexity for special, yet non-negligible many, cases.


\subsection{Organization of the paper}

The paper is organized as follows. In Section~\ref{smallsumofdigits},
we list some results of counting numbers with small sum-of-digits.
In Section~\ref{basicidea}, we present the basic idea and the
algorithm, and prove Theorem~\ref{main}. In
Section~\ref{parameterandlistdecoding}, we prove
Theorem~\ref{listdecoding} and Theorem~\ref{parameter}.
In Section~\ref{artinschreier}, we extend the results
to Artin-Schreier extensions. We conclude
our paper with discussions of open problems.


\section{Numbers with Small Sum-of-digits}\label{smallsumofdigits}

Suppose that the $q$-ary expansion of a positive integer $e$ is
$$e = e_0 + e_1 q + e_2 q^2 + \cdots + e_{n-1} q^{n-1}, $$
where $0\leq e_i \leq q-1$ for all $ 0 \leq i \leq n-1$.
How many nonnegative integers $e$ less than $q^{n}$ satisfy $S_q (e) = w$?
The number equals to the number of nonnegative integral solutions of
$$ \sum_{i=0}^{n-1} e_i = w $$
under the conditions that $0\leq e_i \leq q-1$ for all $ 0 \leq i \leq n-1$.
Denote the number by $N(w,n,q)$.  The generating
function for $N(w,n,q)$ is
$$ ( 1+x+\cdots + x^{q-1} )^n = \sum_i N(i,n,q) x^i. $$ If $w\leq q-1$,
then the conditions $e_i \leq q-1$ can be removed, we have that $
N(w,n,q) = { w + n-1 \choose n-1 }$.  It is easy to see that if $q=2$,
we have that $ N(w,n,2) = {n \choose w} $.  In the later section, we
will need to estimate $N(w,n,q)$, where $w$ is $n$ times a small constant
less than $2$. Since
\begin{eqnarray*} 
&&( 1+x+\cdots + x^{q-1} )^n \\
&=& ( { 1 - x^q \over 1 -x} )^n\\ 
&=& (1 -x^q)^n \sum_{i=0}^\infty { i + n - 1 \choose n-1 } x^i \\ 
&\equiv& (1- n x^q ) \sum_{i=0}^{2q-1} { i + n - 1 \choose n-1 } x^i
\pmod{x^{2q}}\\
&\equiv& \sum_{i=0}^{q-1} { i + n - 1 \choose n-1 } x^i
+ \sum_{i=q}^{2q-1} ( { i + n - 1 \choose n-1 } 
- n { i - q + n-1 \choose n-1 } ) x^i \pmod{x^{2q}}
\end{eqnarray*}

Hence $N(w, n, q) = { w + n - 1 \choose n-1 } 
- n { w - q + n-1 \choose n-1 }$ if $ q\leq  w < 2q$.


\section{The Basic Ideas and the Algorithm}\label{basicidea}

Our basic idea is adopted from the index calculus algorithm.
Let $\F_{q^n}$ be a Kummer extension of $\F_q$, namely, $n|q-1$.
Assume that $q = p^d$ where $p$ is the characteristic.
The field $\F_{q^n}$ is usually given as 
$\F_p[x]/(u(x))$ where $u(x)$ is an irreducible polynomial 
of degree $dn$ over $\F_p$. If $g$ satisfies the condition
in Theorem~\ref{main}, then $x^n - \alpha^n$ must be
an irreducible polynomial over $\F_q$. Denote $\alpha^n $ by $a$.
To implement our algorithm, it is
necessary that we work in
another model of $\F_{q^n} $, namely, $ \F_q[x]/(x^n-a)$. 
Fortunately the isomorphism 
$$\psi: \F_p[y]/(u(y)) \rightarrow \F_{q^n} = \F_q[x]/(x^n-a)$$ can be
efficiently computed. To computer $\psi ( v(y))$, where $v(y)$ is a
polynomial of degree at most $dn-1$ over $\F_p$, all we have to do is
to factor $u(y)$ over $\F_q[x]/(x^n-a)$, and to evaluate $v(y)$ at 
one of the roots. The random algorithm runs in expected time $ O(dn (dn + \log
q^n) (dn \log q^n )^2) $, and 
the deterministic algorithm runs in time $ O(dn (dn + q) (dn
\log q^n )^2) $.  From now on we assume the model $\F_q[x]/(x^n-a)$.

Consider the subgroup generated by $g = \alpha+b$ in $(\F_q[x]/(x^n-a))^*$,
recall that $b \in \F_q^*$ and $\alpha = x \pmod{x^n -a}$.  The
generator $g$ has
order greater than $2^n$ \cite{Cheng04}, and
has a very nice property as follows.  Denote $a^{q-1 \over n}$ by $h$,
we have
$$ g^q = (\alpha+b)^q = \alpha^q + b = a^{q-1 \over n} \alpha + b = h
\alpha + b, $$ and more generally
$$ (\alpha+b)^{q^i} = \alpha^{q^i} + b = h^i \alpha + b. $$ In the
other word, we obtain a set of relations: 
$\log_{\alpha + b} (h^i \alpha + b ) = q^i $
for $0\leq i\leq n-1$.  This
corresponds to the precomputation stage of the index calculus.  The
difference is that, in our case, the stage finishes in polynomial time,
while generally it requires subexponential time. For a
general exponent $e$,
$$ (\alpha+b)^e = (\alpha+b)^{e_0 + e_1 q + \cdots + e_{n-1} q^{n-1}}
= (\alpha+b)^{e_0} (h \alpha + b)^{e_1} \cdots (h^i \alpha+b)^{e_i}
\cdots (h^{n-1} \alpha + b)^{e_{n-1}}. $$ If $  f(\alpha) $ is an
element in $\F_{q^n}$, where $f \in \F_q [x]$ is a polynomial of
degree less than $n$, and $f(\alpha) = (\alpha+b)^e $ and $S_q (e) < n$, 
then due to unique factorization in $\F_q [x]$,
{\em $f(x)$ can be completely split into the product of linear
factors over $\F_q$.  We can read the discrete logarithm from the
factorizations, after  the coefficients
are normalized.}
The algorithm is described as follows.

\begin{algorithm}

Input: $g$, $g^e$ in $\F_{q^n} = \F_q[x]/(x^n-a)$ satisfying the
conditions in Theorem~\ref{main}.

Output: $e$.

\begin{enumerate}
\item \label{computingandsorting} 
Define an  order in $\F_q$ (for example, use the alphabetic order).
Compute and sort the list $(1, h, h^2, h^3, \cdots, h^{n-1})$. 
\item Suppose that $g^e$ is represented by $f(\alpha)$,
where $f\in \F_q [x]$ has degree less than $n$.
Factoring $f(x)$ over $\F_q$, let $f(x) = c (x+d_1)^{e_1} \cdots
   (x+d_k)^{e_k} $ where $c, d_1, \cdots, d_k$ are in $\F_q$.

\item (Normalization) Normalize the coefficients and reorder the
factors of $f(x)$ such that their constant coefficients 
are $b$ and $f(x) = (x+b)^{e_1} \cdots (h_{n-1} x+b)^{e_{n-1}} $,
where $h_i = h^i$;
\item Output $e_0 + e_1 q + \cdots + e_{n-1} q^{n-1}$;
\end{enumerate}

\end{algorithm}

The step~\ref{computingandsorting} takes time $O(n \log^2 q \log n + n
\log n \log q ) = O(n \log n \log^2 q) $.  The most time-consuming
part is to factor a polynomial over $\F_q$ with degree at most
$n$. The random algorithm runs in expected time
$O(n(n+\log q) ( n \log q)^2)$ and
the deterministic algorithm runs in time
$O(n(n+q) ( n \log q)^2) = O( n^3 q \log^2 q)$.
Normalization and reordering
can be done in time $O(n\log n \log q)$, since we have a
sort list of $(1, h, h^2, h^3, \cdots, h^{n-1})$. The total time
complexity is thus in random time $O(n(n+\log q) ( n \log q)^2)$ and 
in deterministic time $O(n^3 q \log^2 q)$.  
This concludes the proof of the main theorem.

\section{The Parameterized Complexity and The Application of
List Decoding}\label{parameterandlistdecoding}

A natural question arises: can we relax the bound on the 
sum-of-digits
and still get a polynomial time algorithm?
Solving the problem under the condition $S_q (e) \leq (q-1)n$
basically  renders the discrete logarithm problems
in $\F_{q^n}$ and any of its subfields easy.
In this section, we consider the case when $S_q (e) \leq 1.32n$.
Suppose that $g^e = f(\alpha)$
where $f(x) \in \F_q [x]$ has degree less than $n$. 
Use the same notations as in the previous section, we have
$$ f(\alpha) = (\alpha+b)^{e_0} (h\alpha + b)^{e_2} \cdots (h^{n-1} \alpha +
b)^{e_{n-1}}. $$ 
Hence there exists a polynomial
$t(x)$ with degree less than $0.32 n$ such that 
$$ f(x) + (x^n - a)t(x)  =
(x+b)^{e_0} (h x + b)^{e_1} \cdots (h^{n-1} x + b)^{e_{n-1}}.$$
If there are at least $0.5657n > \sqrt{0.32n \cdot n}$ 
number of nonzero $e_i$'s,
then the curve $y = t(x)$ will pass at least $0.5657n$ point
in the set
$$\{ (i, - { f(i) \over i^{q-1} - a }) | i \in \{-b, -{b \over h},
\cdots, -{b \over h^{n-1}} \} \}. $$ To find all the polynomials of
degree less than $0.32 n$, which pass at least $0.5657 n$ points in a
given set of $n$ points, is an instance of the list decoding problem.
It turns out that there are only a few of such polynomials, and
they can be found efficiently.

\begin{proposition} (Guruswami-Sudan \cite{GuruswamiSu99} ) 
Given $n$ distinct elements $x_0, x_1, \cdots, x_{n-1} \in \F_q $,
$n$ values $y_0,y_1,\cdots,y_{n-1} \in \F_q $ and a natural number $k$,
there are at most
$O(\sqrt{n^3 k})$ many univariate polynomials $t(x) \in \F_q [x]$ 
of degree at most $k$ such that
$y_i = t(x_i)$ for at least $\sqrt{nk}$ many points.
Moreover, these polynomials can be found in random polynomial time.
\end{proposition}

For each $t(x)$, we use the Cantor-Zassenhaus algorithm
to factor $ f(x) + (x^n -a)* t(x) $.
There must exist a $t(x)$ such that the polynomial $f(x) + (x^n
-a)* t(x)$ can be completely factored into a product of linear factors
in $\{ h^{q^i} x + b | 0\leq i\leq n-1\}$, and $e$ is computed as a
consequence.
In order to prove Theorem~\ref{listdecoding}, it remains
to show:

\begin{lemma}
Define 
$$A_{n,q} =\{(x_1,x_2, \cdots, x_n) \mid x_1 + x_2 + \cdots + x_{n} \leq 1.32n,
x_i \in \Z {\rm\ and\ } 0\leq x_i\leq q-1 \ {\rm\ for\ }\ 1\leq i \leq n. \}  $$ 
and 
$$ B_n = \{(x_1,x_2, \cdots, x_n) \mid |\{ i| x_i = 0\}| \geq 0.5657 n\}. $$
We have
$$ { |A_{n,q} \cap B_n| \over |A_{n,q}|  } < c^{-n} $$
for some constant $c>1$ when $n$ is sufficiently large.
\end{lemma}

\begin{proof}
The cardinality of $A_{n,q}$ is $\sum_{i=0}^{\lfloor 1.32n \rfloor}
N(i,n,q) > { 2.32n \choose n} > 4.883987...^n $.  The cardinality of
$A_{n,q} \cap B_n$ is less than $ \sum_{v= \lceil 0.5657n \rceil }^n
{n \choose v} { 1.32n \choose n-v-1} $. The summands maximize at $v =
0.5657 n $ if $v\geq 0.5657n$.  Hence we have
\begin{eqnarray*}
&&\sum_{v= \lceil 0.5657n \rceil }^n {n \choose v} { \floor{1.32 n}
\choose n-v-1} \\ &<& 0.5657 n {n \choose \ceil{0.5657n}} {\floor{1.32
n} \choose \floor{0.4343 n}}\\ &<& 4.883799...^{n }
\end{eqnarray*}
This proves the lemma with  $ c = 4.883987.../4.883799... > 1$. 
\end{proof}

Now we are ready to prove Theorem~\ref{parameter}.  Any $f(x)$ where
$f(\alpha) = (\alpha+b)^e \in <\alpha+b> \subseteq \F_{q^{q-1}}$ is 
congruent to a product of at most $w = S_q (e)$ linear factors
modulo $x^{q-1} - a $. If $w< q-1$, we have
an algorithm running in time $O(q^4 \log^2 q)$, according to
Theorem~\ref{main}.  So we only need to consider the case when $ w
\geq q-1$.  The general purpose algorithm will run in random
time $f(\log q^{q-1})$, where $f$ is a subexponential function.
Since $\log q^{q-1} \leq w\log w$, this proves
Theorem~\ref{parameter}.


\section{Artin-Schreier Extensions}\label{artinschreier}
Let $p$ be a prime.
The Artin-Schreier extension of a finite field $\F_p$ 
is $\F_{p^p}$. It is easy to show that $x^p - x - a =0$
is an irreducible polynomial in $\F_p$
for any $a \in \F_p^*$. So we may take
$\F_{p^p} = \F_p [x]/(x^p-x-a)$. Let $\alpha = x \pmod{x^p-x-a}$.
For any $b \in \F_p$,
we have
$$ (\alpha + b)^p = \alpha^p + b = \alpha + b + a, $$
and similarly 
$$ (\alpha + b)^{p^i} = \alpha^{p^i} + b = \alpha + b + ia. $$
Hence the results for Kummer extensions can be adopted to
Artin-Schreier extensions.
For the subgroup generated by $\alpha + b$, we have
a polynomial algorithm to solve the discrete logarithm
if the exponent has $p$-ary sum-of-digits less than $p$.
Note that $b$ may be $0$ in this case.

\begin{theorem}
There exists an algorithm to find the integer $e$
given $g$ and $g^e$ in $\F_{p^p}$ in time polynomial in $\log p^p$
under the conditions:
\begin{enumerate}
\item $0\leq e< p^p$, and $S_q (e) \leq p-1$;
\item $g = \alpha+b $ where $\F_p (\alpha) =
\F_{p^p}$, $b\in \F_q $ and $\alpha^p + \alpha \in \F_p^*$.
\end{enumerate}
Moreover, there does not exist an integer $e' \not= e$  
satisfying that $0\leq e'< p^p$, $S_q (e') \leq n$
and $g^{e'} = g^e$.

\end{theorem}

\begin{theorem}
There exists an element $g$ of order greater than $2^p$ in 
$\F_{p^p}^*$, such that the discrete logarithm problem 
with respect to $g$
can be solved in time $O(f(w) (\log p^p)^4)$,
where $f$ is a subexponential function and
$  w $ is the bound of the sum-of-digits of the 
exponent in the $p$-ary expansion.

\end{theorem}

\begin{theorem}
Suppose that $g = \alpha+b $, where $\F_p (\alpha) = \F_{p^p}$, $b\in \F_p $
and $\alpha^p + \alpha \in \F_p^*$.
Suppose $e$ is chosen in random  from the set
$$ \{ 0\leq e < q^n -1 | S_q (e) < 1.32 n \}. $$ There exists an algorithm
given $g$ and $g^e$ in $\F_{p^p}$, to find $e$ in time polynomial in
$\log (p^p)$, with probability greater than $1 - c^{-n}$ for some
constant $c$ greater than $1$.
\end{theorem}

\section{Conclusion Remarks}



A novel idea in the celebrated AKS primality testing algorithm,
is to construct a subgroup of large cardinality through linear elements
in finite fields.
The subsequent improvements \cite{Berrizbeitia02, Cheng03b, Bernstein03}
rely on constructing a single element of large order.
It is speculated that these ideas will be useful in 
attacking the integer factorization problem.
In this paper, we show that they do affect the discrete logarithm
problem in finite fields.
We give an efficient algorithm which computes the 
bounded sum-of-digits discrete logarithm with respect to 
prescribed bases
in Kummer fields. We emphasize that
this is more than a result which deals with only special cases,
as every finite field has extensions of reasonable degrees 
which are Kummer fields.
One of the most interesting problems is to further relax the
restriction on the sum-of-digits of the exponent.
Another important open problem is to prove Conjecture~\ref{orderconjecture}. 
If that conjecture is true, the AKS-style primality
proving can be made compatible or better than
ECPP or the cyclotomic testing in practice.

\paragraph{Acknowledgments}
We thank Professor Pedro Berrizbeitia
for very helpful discussions.

\bibliographystyle{plain}
\bibliography{crypto}

\end{document}